\documentclass{article}
\usepackage{latexsym,amsmath,amssymb,amscd}
\usepackage[dvips]{color}
\begin{document}

\title{On slowly rotating axisymmetric solutions of the
Euler-Poisson equations }
\author{Juhi Jang \thanks{Department of Mathematics, University of Southern California, e-mail: juhijang@usc.edu} and Tetu Makino \thanks{
Professor Emeritus at Yamaguchi University, e-mail: makino@yamaguchi-u.ac.jp}}
\date{\today}
\maketitle

\newtheorem{Lemma}{Lemma}
\newtheorem{Proposition}{Proposition}
\newtheorem{Theorem}{Theorem}
\newtheorem{Definition}{Definition}
\newtheorem{Remark}{Remark}
\newtheorem{Corollary}{Corollary}

\begin{abstract}
We construct stationary axisymmetric solutions of the Euler-Poisson equations, which govern the internal structure of polytropic gaseous stars, with small constant angular velocity when the adiabatic exponent $\gamma$ belongs to $(\frac65,\frac32]$.  The problem is formulated as a nonlinear  integral equation, and is solved by iteration technique. By this method,
not only we get the existence, but also
we clarify properties of the solutions such as the physical vacuum condition and oblateness of the star surface.  \\

{\it Key Words and Phrases.} Axisymmetric solutions. Gaseous stars. 
Euler-Poisson equations. Stellar rotation. Nonlinear Integral equation.
Semilinear elliptic problem. Free boundary.\\

{\it 2010 Mathematics Subject Classification Numbers.} 35Q35, 35Q85, 35J61, 45G10, 76D05.

\end{abstract}

\section{Introduction}

The Euler-Poisson equations which govern the evolution of a gaseous star are
\begin{subequations}
\begin{align}
&\frac{\partial\rho}{\partial t}+\nabla\cdot(\rho\vec{v})=0, \label{B11a}\\
&\rho\Big(\frac{\partial\vec{v}}{\partial t}+
(\vec{v}\cdot\nabla)\vec{v}\Big)+\nabla P=-\rho\nabla\Phi, \\
&\triangle\Phi=4\pi\mathsf{G}\rho. \label{B1c}
\end{align}
\end{subequations}

Here the independent variable is $(t,\vec{x})=(t,x,y,z)\in [0,T)\times\mathbb{R}^3$, $\mathsf{G}>0$ is 
 the gravitational constant, and unknown functions are
the density $\rho=\rho(t,\vec{x})\geq 0$, the velocity field
$\vec{v}=\vec{v}(t,\vec{x})\in \mathbb R^3$, the gravitational potential
$\Phi=\Phi(t,\vec{x})$. We assume that the pressure $P$ is a function of $\rho$ given by the $\gamma$-law:  
\begin{equation}\label{gammalaw}
P=\mathsf{A}\rho^{\gamma},
\end{equation}
where $\gamma\in (1,2)$ is the adiabatic exponent and 
$\mathsf{A}>0$ is a given entropy constant. 
In this article we denote
\begin{equation}
\nu=\frac{1}{\gamma-1},
\end{equation}
whence $1<\nu<+\infty$.

Since we are interested in density distributions with compact support, we specify the
potential satisfying the Poisson equation \eqref{B1c} by the Newtonian potential
\begin{equation}
\Phi(t,\vec{x})=-\mathsf{G}
\int\frac{\rho(t,\vec{x}')}{|\vec{x}-\vec{x}'|}d\mathcal{V}(\vec{x}'), 
\end{equation}
where $d\mathcal{V}$ is the usual volume element $dx'dy'dz'$.\\

The most well-known solutions to the Euler-Poisson equations \eqref{B11a}-\eqref{B1c} are the Lane-Emden functions, which describe non-rotating spherically symmetric stars in equilibrium. The ansatz $\rho=\rho(r)$ and $\vec{v}=\vec{0}$ gives rise to the second order ODE, known as the Lane-Emden equation. The Lane-Emden solutions have been extensively studied by astrophysicists \cite{Ch} and they offer interesting mathematical   problems. The dynamics of spherically symmetric solutions of \eqref{B11a}-\eqref{B1c} with \eqref{gammalaw} have been studied by the authors \cite{MakinoOJM, Jang2014,JangAPDE} near Lane-Emden solutions. 

The Euler-Poisson equations are also used to model rotating stars. In this case, even the existence theory of stationary solutions is not complete. This is because one must deal with more than one variable (PDEs) once the spherical symmetry assumption is dropped, and often it requires more sophisticated arguments.  A key difficulty of obtaining axisymmetric rotating star solutions lies in the fact that the boundary of the domain $\{\vec{x}: \rho(\vec{x})>0\}$ is not prescribed and it is part of the problem: a free boundary problem. 

The aim of this article is to construct a family of compactly supported, stationary axisymmetric solutions of  \eqref{B11a}-\eqref{B1c}  obeying the $\gamma$-law \eqref{gammalaw} with small constant angular velocity which demonstrates slow solid rotations. \\

We look for solutions to  \eqref{B11a}-\eqref{B1c} of the form 
\begin{equation}\label{ansatz}
\rho=\rho(\varpi, z), \;\;\vec{v}=(-\Omega y, \Omega x, 0)^T
\end{equation}
 where the angular velocity $\Omega$ is 
  a constant. 
Here 
$(\varpi, z,\phi)$ denote the cylindrical coordinates: 
\begin{equation}
x=\varpi \cos\phi,\quad
y=\varpi \sin\phi,\quad z=z.
\end{equation}
If the polar coordinates are denoted by $(r,\vartheta, \phi)$: 
$$ x=r\sin\vartheta\cos\phi,\quad y=r\sin\vartheta\sin\phi,\quad
z=r\cos\vartheta,
$$
then 
$$\varpi=r\sin\vartheta=\sqrt{r^2-z^2}.
$$

Instead of the cylindrical coordinates $(\varpi, z)$, we shall use another spherical polar coordinates $(r,\zeta)$  
defined by
\begin{equation}
x=r\sqrt{1-\zeta^2}\cos\phi,\quad y=r\sqrt{1-\zeta^2}\sin\phi,
\quad z=r\zeta,
\end{equation}
and hence 
$$\varpi=r\sqrt{1-\zeta^2},\quad z=r\zeta.
$$
The Euler-Poison equations with the ansatz \eqref{ansatz} are then reduced to
\begin{subequations}
\begin{align}
-\rho(1-\zeta^2)r\Omega^2+\frac{\partial P}{\partial r}&=-\rho
\frac{\partial\Phi}{\partial r},\label{B1a}\\
\rho\zeta r^2\Omega^2+\frac{\partial P}{\partial\zeta}&=
-\rho\frac{\partial\Phi}{\partial\zeta}, \label{B1b}\\
\triangle\Phi&=4\pi\mathsf{G}\rho,
\end{align}
\end{subequations}
where
\begin{equation}
\triangle\Phi=
\frac{1}{r^2}\frac{\partial}{\partial r}r^2\frac{\partial\Phi}{\partial r}
+
\frac{1}{r^2}\frac{\partial}{\partial\zeta}(1-\zeta^2)\frac{\partial\Phi}{\partial\zeta}.
\end{equation}
The Newton potential reads
\begin{equation}
\Phi(r,\zeta)=-\mathsf{G}
\int_{-1}^1\int_0^{\infty}
K(r,\zeta,r',\zeta')\rho(r',\zeta')r'^2dr'd\zeta',
\end{equation}
where
\begin{equation}
K(r,\zeta,r',\zeta')=\int_0^{2\pi}
\frac{d\beta}{\sqrt{r^2+r'^2-2rr'
(\sqrt{1-\zeta^2}\sqrt{1-\zeta'^2}\cos\beta+\zeta\zeta')}}.
\end{equation}

\

It is convenient to introduce the enthalpy variable $u$ defined by
\begin{equation}\label{enthalpy}
u=\int_0^{\rho}\frac{dP}{\rho}=\frac{\mathsf{A}\gamma}{\gamma-1}\rho^{\gamma-1}.
\end{equation}
On the region where $\rho>0$, the system \eqref{B1a} \eqref{B1b}
is then reduced to
\begin{align*}
-(1-\zeta^2)r\Omega^2+\frac{\partial u}{\partial r}&=
-\frac{\partial\Phi}{\partial r} ,\\
\zeta r^2\Omega^2+\frac{\partial u}{\partial \zeta}&=
-\frac{\partial\Phi}{\partial\zeta},
\end{align*}
which is equivalent to
$$\Phi+u=\frac{\Omega^2}{2}(1-\zeta^2)r^2+\mbox{C},
$$
where $C$ is a constant. 
By a change of scales, without loss of generality, we may assume
$$ 4\pi\mathsf{G}\rho=u_{\sharp}^{\nu}, $$
where we denote 
$$u_{\sharp}(r,\zeta):=\max\{0, u(r,\zeta)\}$$
when $u$ is extended negatively outside the support of $\rho$.
We set $\Omega^2=\varepsilon/2$. 

We now formulate the problem:\\

{\bf (P) \ Find 
a bounded domain $\mathfrak{D}$ including the origin $\vec{0}$ and
a function 
$u=u(r,\zeta)$ which is defined and positive on
the 
domain $\mathfrak{D}$, continuous on the closure of $\mathfrak{D}$, 
vanishing on the boundary $\partial\mathfrak{D}$ of the domain $\mathfrak{D}$, and
satisfies
\begin{equation}\label{integralE}
\Phi(r,\zeta)+u(r,\zeta)=
\frac{\varepsilon}{4}(1-\zeta^2)r^2+u_c+
\Phi(\vec{0})
\end{equation}
in the domain $\mathfrak{D}$ with the potential given by
\begin{equation}
\Phi(r,\zeta)=-\frac{1}{4\pi}
\int_{-1}^1\int_0^{\infty}
K(r,\zeta,r',\zeta')u_{\sharp}(r',\zeta')^{\nu}r'^2dr'd\zeta',
\end{equation}
where $u_{\sharp}=0$ outside $\mathfrak{D}$.} \\

Here $u_c=u(\vec{0})>0$ is the central enthalpy value 
and we may assume $u_c=1$ without loss of generality. 
Note that
\begin{align*}
\Phi(\vec{0})&=-\frac{1}{4\pi}\int_{-1}^1\int_0^{\infty}K(0,0,r',\zeta')u_{\sharp}(r',\zeta')^{\nu}
r'^2dr'd\zeta' \\
&=-\frac{1}{2}\int_{-1}^1\int_0^{\infty}u_{\sharp}(r',\zeta')^{\nu}r'dr'd\zeta'.
\end{align*}
We see that $u$ satisfies 
\begin{equation}
-\triangle u=u^{\nu}-\varepsilon \ \text{ in } \  \mathfrak{D}. 
\end{equation}\\

The first attempt to construct axisymmetric rotating stationary solutions for the Euler-Poisson equations \eqref{B11a}-\eqref{B1c} with \eqref{gammalaw}
was made by E. A. Milne \cite{Milne} in 1923 for the case of
$\gamma=4/3$ ($\nu=3$), and by
S. Chandrasekhar  \cite{Chandrasekhar1933} in 1933 for general $\gamma$ based on the perturbation method. 
Assuming that 
$\varepsilon$ is a sufficiently small constant, 
S. Chandrasekhar tried to find stationary solutions of the form
$$u=\sum_{n=0}^{\infty} u_n(r,\zeta)\varepsilon^n $$
with 
$u_0(r,\zeta)=u_c\theta(r/\alpha  ;\  \nu),$
$u_c, \alpha$ being positive constants
and $\theta$ being the Lane-Emden function. Some arguments, however, are not mathematically rigorous and contain logical gaps, although it is certainly a pioneering work.

On the other hand, a mathematically rigorous treatment of the problem was initiated by
J. F. G. Auchmuty and R. Beals \cite{AuchmutyB} in 1971.
They established the method to construct solutions through
the variational problem to minimize the energy under the given total mass and
the total angular momentum, which should be conserved along the evolutions.
Along this line, many mathematically rigorous studies appeared \cite{Ach1991, CaF,ChLi,FrTur1,FrTur3,LuoS0, LuoS}. In \cite{Li}, Y. Y. Li  constructed uniformly rotating solutions with small constant angular velocity based on the variational formulation. We remark that the results obtained by variational approaches require specific assumptions on the equation of state, and at least the exact $\gamma$-law \eqref{gammalaw} with $\gamma <4/3$ does not seem to be covered by those results.  

We will study the problem by a non-variational approach,
which is a natural justification of the method adopted by astrophysicists
after S. Chandrasekhar.

\section{Main results}\label{Sec:2}

\subsection{Lane-Emden function}

We first assume \\ 

{\bf (AS0):} \hspace{10mm} $1<\nu<5$, that is, $\frac{6}{5}<\gamma<2$.\\

Let $\theta=\theta(r)=\theta(r;\nu)$ denote the Lane-Emden function of index $\nu$. 
That is, it is the unique solution of 
\begin{equation}\label{LE}
\frac{d^2u}{dr^2}+\frac{2}{r}\frac{du}{dr}+u^{\nu}=0,\quad 
u=1+O(r^2)\quad\mbox{as}\quad r\rightarrow 0.
\end{equation}

It is well-known \cite{Ch} that under {\bf (AS0)} there exists a finite $\xi_1(\nu) (>\pi)
$ such that $\theta(r;\nu)>0, \frac{d}{dr}\theta(r;\nu)<0$ for $0<r<\xi_1(\nu)$,
and $\theta(\xi_1(\nu);\nu)=0$. \\

We consider the extension of $\theta(r;\nu)$, denoted by $\theta \in C^2([0,+\infty))$, such that 
\begin{equation}\label{extension}
\theta(r)=-\mu_1(\nu)\Big(\frac{1}{\xi_1(\nu)}-\frac{1}{r}\Big) \ \text{ for } \ r >\xi_1(\nu)
\end{equation}
where 
$$\mu_1(\nu):=\int_0^{\xi_1(\nu)}\theta(r;\nu)^{\nu}r^2dr=
-r^2\frac{d\theta(r;\nu)}{dr}\Big|_{r=\xi_1(\nu)}$$
which is proportional to the total mass of the star.  
Then $\triangle\theta =0$ and $\theta<0$ for 
$r>\xi_1(\nu)$.  Note that $\theta(r)$ solves 
\[
\frac{d^2u}{dr^2}+\frac{2}{r}\frac{du}{dr}+u_\sharp^{\nu}=0
\]
for all $r<+\infty$. The harmonic extension by \eqref{extension} that gives negative values for $\theta$ for $r>\xi_1(\nu)$ will play an important role in our analysis. In what follows, $\theta$ denotes the above extension.\\

Here are some notations used throughout the paper. 
\begin{itemize}
\item We denote $
u_{\sharp}(r,\zeta):=\max\{0, u(r,\zeta)\}.
$
\item 
$
 a \wedge b:=\min\{a,b\}\;\text{ for }\;a,b\in \mathbb R. 
$
\item We denote 
$
[X]_1 := \sum_{k\geq 1}a_kX^k 
$ 
with positive radius of convergence. For example, we note that $\theta(r)=1+[r^2]_1$ as $r\rightarrow 0$. 
\item We denote  
$
\nabla =\left(\frac{\partial}{\partial x}, \frac{\partial}{\partial y}, \frac{\partial}{\partial z}\right)^T. 
$
\end{itemize}

\subsection{Main results} 

Let us fix $R_0\geq 2\xi_1(\nu)$.  By concentrating on solutions which are symmetric with respect to
the plane $\{ z=0\}$, we shall seek solutions in the 
 function space
\begin{align}
\mathfrak{E}:=&\{ u\in C([0,R_0]\times[-1,1]) : \nonumber \\
&\;\;u(0,0)=u(0,\zeta)\;\; 
\mbox{and}\;\; u(r,\zeta)=u(r,-\zeta)\;\;\mbox{for}\;\; \forall \zeta\in[-1,1]\}
\end{align}
endowed with the norm
\begin{equation}
\|u\|_{\mathfrak{E}}:=\sup |u(r,\zeta)|.
\end{equation} 

We now state our main results. 

\begin{Theorem}\label{Theorem1} Assume $2\leq \nu<5$ or $\frac65<\gamma\leq \frac32$. Then 
there exists $\varepsilon_0>0$ such that for each $0<\varepsilon\leq \varepsilon_0$, there exists $\Theta=\Theta(r,\zeta;\nu,\varepsilon) \in\mathfrak{E}$  
satisfying the nonlinear integral equation
\begin{equation}
u(r,\zeta)=\varepsilon\mathfrak{g}(r,\zeta)+\mathcal{G}(u)(r,\zeta),\label{6}
\end{equation}
for each $(r,\zeta)\in [0,R_0]\times[-1,1]$. Here 
\begin{equation}
\mathfrak{g}(r,\zeta):=\frac{1}{4}(1-\zeta^2)r^2
\end{equation}
and
\begin{align}
\mathcal{G}(u)(r,\zeta):=&
1+\frac{1}{4\pi}\int_{-1}^1\int_0^{R_0}
K(r,\zeta,r',\zeta')u_{\sharp}(r',\zeta')^{\nu}r'^2dr'd\zeta' 
+ \nonumber \\
&-\frac{1}{2}\int_{-1}^1\int_0^{R_0}
u_{\sharp}(r',\zeta')^{\nu}r'dr'd\zeta'.\label{8}
\end{align}
Moreover, it holds that $\| \Theta-\theta\|_{\mathfrak{E}}\leq C \varepsilon$ for some constant $C>0$. 
\end{Theorem}

The solution $\Theta$ will be referred to as the distorted
Lane-Emden function of index $\nu$ with the squared angular
velocity parameter $\varepsilon$. 

\begin{Theorem}\label{Theorem2}  The distorted
Lane-Emden function $\Theta$ in Theorem \ref{Theorem1} satisfies the following properties. 
\begin{enumerate}
\item[(1)] $\Theta$ is continuously differentiable
as a function of  $\vec{x}$, where $$\vec{x}=(x,y,z)=(r\sqrt{1-\zeta^2}\cos\phi, r\sqrt{1-\zeta^2}\sin\phi, r\zeta), $$
namely $\nabla \Theta \in \mathfrak{E}$ with $\|\nabla \Theta -\nabla \theta\|_\mathfrak{E}\leq C \varepsilon$ for some constant $C>0$. Moreover, $\Theta$ is twice differentiable and the second derivatives are H\"{o}lder continuous
as a function of $\vec{x}$ . 
\item[(2)] 
There exists a unique $\Xi_1(\zeta)=\Xi_1(\zeta;\nu,
\varepsilon)$ in  $(0,2\xi_1)$ for each $\zeta\in [-1,1]$ such that
\begin{align}
&\Theta(r,\zeta;\nu,\varepsilon)>0 \quad\mbox{for}\quad
0\leq r <\Xi_1(\zeta), \nonumber \\
&\Theta(\Xi_1(\zeta),\zeta;\nu,\varepsilon)=0, \nonumber \\
&\Theta(r,\zeta;\nu,\varepsilon)<0\quad\mbox{for}\quad
\Xi_1(\zeta)<r\leq 2\xi_1 \label{B3}
\end{align}
for $0<\varepsilon\leq\epsilon_0$.  Moreover, $\Xi_1(\zeta)$ is continuously differentiable in $\zeta\in(-1,1)$. 
\item[(3)] $\Theta$ satisfies the physical vacuum boundary condition: 
\[
-\infty < \frac{\partial \Theta}{\partial \vec N}=\nabla \Theta \cdot \vec{N}  <0
\]
where $\vec N$ denotes the unit outer normal vector along $\partial\mathfrak{D}$. 
\item[(4)] The boundary $\partial\mathfrak{D}= \{(r,\zeta)|\Theta(r,\zeta;\nu,\varepsilon)=0\}$ is a $C^1$-submanifold of $\mathbb R^3$. 
\end{enumerate}
\end{Theorem}

Theorem \ref{Theorem1} and Theorem \ref{Theorem2} provide a solution to the problem {\bf (P)} proposed in the introduction for $\frac65<\gamma\leq \frac32$. 
To the best of our knowledge, Theorem \ref{Theorem1} is the first rigorous result on the existence of axisymmetric solutions for the Euler-Poisson equations with the constant angular velocity for the exact $\gamma$-pressure law at least for $\gamma<\frac43$. Moreover, they give a rigorous justification of astrophysicists' early works \cite{Milne,Chandrasekhar1933}. 

We prove  Theorem \ref{Theorem1} by finding nonlinear perturbations around the Lane-Emden functions through  the contraction mapping principle on a subspace of $\mathfrak E$. Because the boundary of the domain of the solutions is not known a priori and it is no longer given by $\xi_1(\nu)$, the solutions are sought in a larger domain by using the negatively extended Lane-Emden function $\theta$. A key for the iteration technique to work is the invertibility of the linearization of $u- \mathcal G(u)$ at $u=\theta$. This requires a quantitative result on the Lane-Emden function $\theta$, which was available only by numerics prior to this paper. We provide a rigorous proof that works for $2\leq \nu<5$ (Theorem \ref{Thm1} in Section \ref{Sec:4}). Theorem \ref{Theorem2} follows mostly from direct computations by using the information and the estimates obtained in order to prove Theorem \ref{Theorem1}. 

The paper proceeds as follows. Section \ref{Sec:3} and Section \ref{Sec:4} are devoted to the study of the Fr\'{e}chet derivative $D\mathcal G$ of $\mathcal{G}$. In particular, in Section \ref{Sec:4},  we show that $1-D\mathcal G$ has a bounded inverse, and also the results of astrophysicists' papers \cite{Milne,Chandrasekhar1933,Kovetz} are discussed. Section \ref{Sec:5} contains the construction of the distorted Lane-Emden function solving the integral equation \eqref{6} and Section 6 contains the proof of Theorem \ref{Theorem2}. 
 In Section \ref{sec:justify}, we will give a mathematical account of Chandrasekhar's approximation given in \cite{Chandrasekhar1933}  on the oblateness of the star surface. In Section \ref{Sec:dis}, we briefly discuss the restriction on $\nu$.

\section{Existence of the Fr\'{e}chet derivative $D\mathcal G$ of
$\mathcal{G}$}\label{Sec:3}

In this section, we will show that the Fr\'{e}chet derivative of
$\mathcal{G}$ exists in $\mathfrak{E}$. To this end, we first present elementary results on the functions in $\mathfrak E$.

\begin{Lemma}\label{Prop1} 
There exists $\delta_0>0$ such that if $u\in\mathfrak{E}$ satisfies $\|u-\theta\|_{\mathfrak{E}}\leqq\delta_0$,
then $u_{\sharp}(r,\zeta)=0$ for $r\geq \frac{3}{4}R_0$.
\end{Lemma}

\noindent{\it Proof.} We first observe that $$\theta(r)\leq \theta( \frac{3}{4}R_0)\leq \theta(\frac32\xi_1(\nu))= -\frac{\mu_1(\nu)}{3\xi_1(\nu)}\;\;\text{for}\;\; r\geq \frac{3}{4}R_0.$$ The result follows by taking $\delta_0=\frac{\mu_1(\nu)}{3\xi_1(\nu)}$.  $\square$  \\

Let us fix $\Lambda_0\geq 2\max\{1, \mu_1(\nu)/\xi_1(\nu)\}$. We then have

\begin{Lemma}\label{Prop2} There exists $\delta_0>0$ such that 
if $u\in\mathfrak{E}$ satisfies
$\|u-\theta\|_{\mathfrak{E}}\leq \delta_0$, then
$\|u\|_{\mathfrak{E}}\leq\Lambda_0$.
\end{Lemma}

\noindent{\it Proof.} Since $|u(r,\zeta)|\leq |\theta(r)|+ |u(r,\zeta)-\theta(r)|$ and since $|\theta(r)|\leq \max\{1, \frac{\mu_1(\nu)}{\xi_1(\nu)}\}$, the result follows by taking $\delta_0= \max\{1, \frac{\mu_1(\nu)}{\xi_1(\nu)}\}$. 
$\square$  \\

From now let us fix $\delta_0>0$ so that the conclusions of Lemma \ref{Prop1} and Lemma \ref{Prop2} are valid.

\begin{Proposition}\label{Proposition1} Let $\mathcal{G}(u)$ be given as in \eqref{8}. The following holds. 
\begin{enumerate}
\item[(1)] If $u\in\mathfrak{E}$ satisfies
$\|u-\theta\|_{\mathfrak{E}}\leq \delta_0$, then $\mathcal G(u)\in \mathfrak E$. Moreover, we have 
\begin{equation}\label{est1}
\|\mathcal{G}(u)\|_{\mathfrak{E}}\leq C\|u\|_{\mathfrak{E}}^{\nu}
\end{equation}
for some constant $C>0$. 
 \item[(2)] Let $u\in \mathfrak E$ with $\|u-\theta\|_\mathfrak E \leq \delta_0/2$ be given. Then there exists a bonded linear operator $D \mathcal G(u)$ of $\mathcal G$ given by 
 \begin{align}
( D\mathcal{G}(u)h)(r,\zeta)=&
\frac{\nu}{4\pi}\int_{-1}^1\int_0^{R_0}
K(r,\zeta,r',\zeta')u_{\sharp}(r',\zeta')^{\nu-1}
h(r',\zeta')r'^2dr'd\zeta' \nonumber \\
&-\frac{\nu}{2}
\int_{-1}^1\int_0^{R_0}u_{\sharp}(r',\zeta')^{\nu-1}h(r',\zeta')r'dr'd\zeta',  \label{defDg}
\end{align}
with the estimate
\begin{equation}
\|D\mathcal{G}(u)h\|_{\mathfrak{E}}\leq C
\|u\|_{\mathfrak{E}}^{\nu-1}\|h\|_{\mathfrak{E}}, \label{est2}
\end{equation}
for some constant $C>0$. 
 \end{enumerate}
\end{Proposition}

Before we give a proof, we introduce further notations. Let us denote 
\begin{align*}
\mathfrak{C}(R):= \{ &\rho \in C([0,+\infty)\times[-1,1]) : \rho(0,\zeta)=\rho(0,0),\\
&\;\rho (r,\zeta)=\rho (r,-\zeta)\;\;\mbox{for}\;\; \forall \zeta\in[-1,1], 
\;\;\text{and} \;\;\rho(r,\zeta)=0\;\; \text{for}\;\; r\geq R \}
\end{align*}
and we use $\mathfrak{C}_0$ to denote the union of $\mathfrak{C}(R)$, $0<R<+\infty$. 

For $\rho \in \mathfrak{C}_0$, 
we then define the function $\mathcal{K}\rho$ by
\begin{equation}\label{K}
\mathcal{K}\rho(r,\zeta):=\frac{1}{4\pi}
\int_{-1}^1\int_0^{+\infty}K(r,\zeta, r',\zeta')\rho(r',\zeta')r'^2dr'd\zeta'.
\end{equation}
This is nothing but the Newtonian potential, that is, we have
$$-\triangle U=\rho\quad\mbox{for}\quad U=\mathcal{K}\rho,$$
provided that $\rho$ is H\"{o}lder continuous. 

For any $\rho\in \mathfrak{C}_0$, $\mathcal{K}\rho$ is continuous on
$[0,+\infty)\times [-1,1]$ and satisfies $\mathcal{K}\rho(0,\zeta)=\mathcal{K}\rho(0,0)$ for $\forall \zeta\in[-1,1]$,
\begin{equation}\label{estK}
|\mathcal{K}\rho(r,\zeta)|\leq \frac{C}{1+r}\|\rho\|_{L^{\infty}},
\end{equation}
where we can take the constant $C$
depending on $R$ for $\rho\in \mathfrak{C}(R)$. \\

\noindent{\it Proof of of Proposition \ref{Proposition1} (1).}  Using $\mathcal K$ notation, we can rewrite \eqref{8} as
\begin{equation}
\mathcal{G}(u)=1+\mathcal{K}u_{\sharp}^\nu-(\mathcal{K}u_{\sharp}^\nu)(0,0),
\end{equation}
provided that $\|u-\theta\|_{\mathfrak{E}}\leq\delta_0$,
whence $u_{\sharp}(r,\zeta)=0$ for $r\geq 3R_0/4$. 

Therefore from \eqref{estK} we deduce the estimate \eqref{est1} and $\mathcal{G}(u)\in\mathfrak{E}$ for 
$u\in \mathfrak{E}$ with $\|u-\theta\|_{\mathfrak{E}}\leq\delta_0$. $\square$\\

In order to prove the second assertion of Proposition \ref{Proposition1}, we will need the following lemma. 
 
\begin{Lemma}\label{Prop3}
Let $|u|, |u+h|\leq\Lambda_0$, and put
\begin{equation}
(u+h)_{\sharp}^{\nu}-u_{\sharp}^{\nu}=
\nu u_{\sharp}^{\nu-1}h+\mathfrak{R}(u,h).
\end{equation}
Then there is a constant C (depending on $\Lambda_0$) such that
\begin{equation}
|\mathfrak{R}(u,h)|\leq C|h|^{\nu\wedge 2}.
\end{equation}
\end{Lemma}

\noindent{\it Proof of of Proposition \ref{Proposition1} (2).} As a direct consequence of Lemma \ref{Prop3}, we have the Fr\'{e}chet derivative $D\mathcal{G}(u)$ of $\mathcal{G}$ at
$u\in\mathfrak{E}$ such that $\|u-\theta\|_{\mathfrak{E}}\leq\delta_0/2$ given in \eqref{defDg}. 
In other notations, we may write 
\begin{equation}
D\mathcal{G}(u)h=\mathcal{K}(\nu u_{\sharp}^{\nu-1}h)-
(\mathcal{K}(\nu u_{\sharp}^{\nu-1}h))(0,0),\label{FD}
\end{equation}
where $\mathcal K$ is defined in \eqref{K}. From \eqref{estK}, the desired result follows. $\square$\\ 

It now remains to prove Lemma \ref{Prop3}. \\

\noindent{\it Proof of Lemma \ref{Prop3}.} It is easy to verify
$$|(1+x)^{\nu-1}-1|\lesssim
\begin{cases}
|x|\quad\mbox{for}\quad |x|\leq 1 \\
x^{\nu-1}\quad\mbox{for}\quad 1\leq x.
\end{cases}
$$
Then it follows that
$$|(1+x)^{\nu}
-1-
\nu x|\lesssim
\begin{cases}
|x|^2\quad\mbox{for}\quad |x|\leq 1 \\
x^{\nu}\quad\mbox{for}\quad 1\leq x.
\end{cases}
$$

 Now we consider
$$\Delta\rho:=(u+h)_{\sharp}^{\nu}-u_{\sharp}^{\nu}.
$$

Case-(00): Suppose $u>0$ and $u+h>0$. Then
\begin{align*}
\Delta\rho&=(u+h)^{\nu}-u^{\nu}=u^{\nu}\Big(\Big(1+
\frac{h}{u}\Big)^{\nu}-1\Big)\Big)= \\
&=\begin{cases}
u^{\nu}\Big(
\nu \frac{h}{u}+O\Big(\Big|\frac{h}{u}\Big|^2\Big)\Big) 
\quad\mbox{for}\quad |h|\leq u \\
u^{\nu}\Big(
\nu\frac{h}{u}+O\Big(\Big|\frac{h}{u}\Big|^{\nu}\Big)\Big)
\quad\mbox{for}\quad u\leq h
\end{cases}\\
&=\begin{cases}
\nu u^{\nu-1}h+O(u^{\nu-2}|h|^2)\quad\mbox{for}\quad |h|\leq u \\
\nu u^{\nu-1}h+O(|h|^{\nu})\quad\mbox{for}\quad u\leq h
\end{cases}\\
&=\nu u^{\nu-1}h+O(|h|^{\nu\wedge 2}),
\end{align*}
where we have used that 
$|u|^{\nu-2}~\leq\Lambda_0^{\nu-2}$
when $\nu \geq 2$.

Case-(01): Suppose $u>0$ but $u+h\leq 0$. Then
$\Delta\rho=-u^{\nu}$. But $0<u\leq -h$ implies
$|u|\leq |h|$ and
$$|\Delta\rho|\leq |h|^{\nu}\lesssim |h|^{\nu\wedge 2},$$
provided that $|h|\leq|u+h|+|u|\leq2\Lambda_0$. On the other hand,
$$|\nu u_{\sharp}^{\nu-1}h|\leq \nu u^{\nu-1}|h|\leq\nu |h|^{\nu}.$$

Case-(10): Suppose $u\leq 0$ but $u+h>0$. Then
$\Delta\rho=(u+h)^{\nu}, 0<u+h\leq h$,
and $u_{\sharp}=0$.

Case-(11): Suppose $u\leq 0$ and $u+h\leq 0$. Then
$\Delta\rho=0$ and $u_{\sharp}=0$, no problem. $\square$

\section{Eigenvalue problem for 
$D\mathcal{G}$}\label{Sec:4}

Since
$$\frac{\partial}{\partial x}\int
\frac{\rho(x')}{|\vec{x}-\vec{x}'|}d\mathcal{V}(\vec{x}')=-
\int\frac{x-x'}{|\vec{x}-\vec{x}'|^3}\rho(\vec{x'})d\mathcal{V}(
\vec{x}')$$
for $\rho \in \mathfrak E_0$, we have 
$$\|\mathcal{K}\rho\|_{C^1}\leq C\|\rho\|_{L^{\infty}},$$
with a constant $C$ depending on $R$, provided that $\rho \in \mathfrak{C}(R)$, where
\begin{align}
\|f\|_{C^1}:&=\|f\|_{L^{\infty}}+\|\nabla f\|_{L^{\infty}} \nonumber \\
&=\|f\|_{\mathfrak{E}}+
\Big\|\sqrt{\Big(\frac{\partial f}{\partial r}\Big)^2+\frac{1-\zeta^2}{r^2}\Big(\frac{\partial f}{\partial \zeta}\Big)^2}\Big\|_{\mathfrak{E}}.
\end{align}

Hence we see
\begin{equation}
\|D\mathcal{G}(u)h\|_{C^1}\leq C\|u\|_{\mathfrak{E}}^{\nu-1}
\|h\|_{\mathfrak{E}}.\label{18}
\end{equation}
Therefore, by Ascoli's theorem, we know that 
$D\mathcal{G}(u)$ is a compact operator, and $1-D\mathcal{G}(u)$
has a bounded linear inverse if $1$ is not an eigenvalue of
$D\mathcal{G}(u)$. (See \cite[Theorem 6.26]{Kato}.)
We want to claim that it is the case for
$u=\theta$, that is, we want to prove

\begin{Lemma}\label{Lemma1}
If $h\in\mathfrak{E}$ satisfies $D\mathcal{G}(\theta)h=h$,
then $h=0$.
\end{Lemma}

We can reduce the proof of Lemma \ref{Lemma1} to the following

\begin{Lemma}[Milne-von Zeipel-Chandrasekhar]\label{Lem1}
Let $j=1,2,\cdots$and $y=\psi(r)$ be the solution of the equation
$$\Big[-\frac{1}{r^2}\frac{d}{dr}r^2\frac{d}{dr}+
\frac{j(j+1)}{r^2}\Big]y=\nu\theta_{\sharp}(r)^{\nu-1}y
\leqno(E_j)
$$
such that $y \sim r^j$ as $r\rightarrow 0$. Then
$$\frac{j+1}{r}y+\frac{dy}{dr}=0\quad\mbox{at}\quad r=\xi_1(\nu) \leqno(D_j) $$
does not hold for $y=\psi(r)$.
\end{Lemma}

Having little knowledge, we have not yet been able to find literatures which describe rigorous proof of this Lemma \ref{Lem1}. Actually, E. A. Milne (1923) used this fact tacitly
in \cite[p. 134]{Milne} for $\nu=3$, and S. Chandrasekhar (1933) did so in
\cite[p. 395]{Chandrasekhar1933} for general $\nu$ than $3$.
H. von Zeipel (1924) claimed this fact explicitly, but wrote
``The proof is omitted here" in \cite[p. 693]{Zeipel}. 
It seems that many astrophysicists believed 
this without a rigorous proof. 

For the time being, {\it assuming that Lemma \ref{Lem1} is true}, we are going to give a proof of Lemma \ref{Lemma1}.\\

From $D\mathcal{G}(\theta)h=h$, we deduce that 
$$
-\triangle h=\nu\theta_{\sharp}^{\nu-1}h, \qquad
h(0,0)=h(0,\zeta)=0 \quad\mbox{for}\quad \forall \zeta\in [-1,1],
$$
and we may assume that $h \in C^2(\mathbb{R}^3)$ and
$-\triangle h=0$ as $|\vec{x}|>\xi_1(\nu)$. Hence, the following identity should also hold: 
\begin{equation}\label{projection}
\int_{-1}^1 -\triangle h P_j d\zeta=\int_{-1}^1 \nu\theta_{\sharp}^{\nu-1}h P_j d\zeta
\end{equation}
for each Legendre's polynomial, $j\geq 0$. By letting 
$$h_j(r)
:=\Big(j+\frac{1}{2}\Big)\int_{-1}^1h(r,\zeta)P_j(\zeta)d\zeta,
$$
the identity \eqref{projection} leads to the equation $(E_j)$ in Lemma \ref{Lem1} for $y=h_j(r)$ for each $j\geq 0$. Also, since $h(0,\zeta)=0$ for all $\zeta\in [-1,1]$, we have $h_j(0)=0$. 

Since $(E_0)$ has a fundamental system of solutions
$y=\psi_{01}(r), \psi_{02}(r)$ such that 
$\psi_{01}(r)=1+[r^2]_1, \psi_{02}(r)=\frac{1}{r}(1+[r^2]_1)$ as $r\rightarrow
0$, $h_0(0)=0$ implies $h_0\equiv 0$. 

Note that since $h(r,\zeta)=h(r,-\zeta)$, $h_{2k-1}\equiv 0$ for $k\in \mathbb N$. Now let $j$ be an even positive integer. The equation $(E_j)$ has a fundamental system of
solutions $y=\psi_{j1}(r), \psi_{j2}(r)$ such that
$$\psi_{j1}(r)=r^j(1+[r^2]_1),\quad
\psi_{j2}(r)=r^{-j-1}(1+[r^2]_1)$$
as $r\rightarrow 0$. Therefore there exists a constant $C_j$
such that $h_j(r)=C_j\psi_{j1}(r)$.

On the other hand, since $h=O(1)$ as $r\rightarrow\infty$
and $\triangle h=0$ on $r>\xi_1(\nu)$, there are constants $A_j$
such that $h_j(r)=A_jr^{-j-1}$ as $r>\xi_1(\nu)$.
Since $h_j\in C^1$ across $r=\xi_1(\nu)$, we know
\begin{subequations}
\begin{align}
C_j y&=A_jr^{-j-1} \label{17a} \\
C_j\frac{dy}{dr}&=-(j+1)A_j r^{-j-2} \label{17b} 
\end{align}
\end{subequations}
should hold for $y=\psi_{j1}$ at $r=\xi_1(\nu)$.
But by Lemma \ref{Lem1} we know the determinant of the coefficient matrix of the simultaneous equations \eqref{17a} \eqref{17b} for $(C_j, A_j)$
$$\Big(\frac{j+1}{r}y+\frac{dy}{dr}\Big)r^{-j-1}\quad
\mbox{for}\quad y=\psi_{j1}\quad\mbox{as}\quad r=\xi_1(\nu)$$ does not vanish and implies
$C_j=A_j=0$, that is, $h_j\equiv 0$. 

Since $h_j\equiv 0$ for all $j\geq 0$, by Stone-Weierstrass theorem, we conclude that $h\equiv 0$. $\square$

\begin{Remark}\label{Rem1}
A. Kovetz \cite{Kovetz} in 1968 gave a half logically rigorous half experimentally plausible proof of the affair of Lemma \ref{Lem1}.

Let $j\in \mathbb N$ be given and let the solution
$y=\psi(r)$ of $(E_j)$ be such that $y \sim r^j (r \rightarrow 0)$
satisfy $(D_j)$. Multiplying $(E_j)$ by $y$ and integrating it on the
interval $[0,\xi_1(\nu)]$, we find
\begin{equation}
(j+1)ry^2|_{r=\xi_1(\nu)}+
\int_0^{\xi_1(\nu)}\Big(\frac{dy}{dr}\Big)^2r^2dr+
j(j+1)\int_0^{\xi_1(\nu)}y^2dr=
\int_0^{\xi_1(\nu)}\nu\theta^{\nu-1}y^2r^2dr, \label{4.1}
\end{equation}
by integration by parts and $(D_j)$. Therefore, if 
\begin{equation}
\nu\theta(r;\nu)^{\nu-1}r^2 <j(j+1)
\quad\mbox{for}\quad 0\leq r\leq\xi_1(\nu),\label{4.2}
\end{equation}
then \eqref{4.1} cannot hold, that is, $(D_j)$ would lead a contradiction.
In this view, A. Kovetz calculated
\begin{equation}
\bar{m}(\nu):=\sup\{ \nu\theta(r;\nu)^{\nu-1}r^2 | 0\leq r\leq \xi_1(\nu)\}
\label{4.3} 
\end{equation}
for some index $\nu$. The `TABLE 1' of \cite[p.1002]{Kovetz} reads
$$
\begin{array}{r||r|r|r|r|r|r|}
\nu & 1&2&2.5&3&4&5\\ \hline
\bar{m}(\nu) & 9.9&1.8&2.3&4.0&3.8&2.8
\end{array}
$$
According this numerical result, it seems that \eqref{4.2} holds for
$j\geq 2 (\Leftrightarrow j(j+1) \geq 6)$, provided that $\nu \geq 2$,
whence $(D_j)$ cannot hold as claimed by Lemma \ref{Lem1}.
But $\bar{m}(\nu)$ can exceed $2$($=j(j+1)$ for $j=1$), and, if we adopt this logic only, we cannot be sure that $(D_1)$ does not hold.
However, we need $(D_j)$ only even $j=2,4,\cdots.$
So, the numerical experiment by A. Kovetz provides an evidence of Lemma \ref{Lemma1} for $\nu\geq2$. 
\end{Remark}

\begin{Remark}
The condition \eqref{4.2}, which is plausible for $j\geq 2, \nu \geq 2$ thanks to the numerical 
result by \cite{Kovetz}, deduces another property of the solution
$y=\psi(r)$ of $(E_j)$ such that $y =r^j(1+[r^2]_1)$ as $ r\rightarrow 0$. That is, 
we find that $\psi(r)$ is monotone increasing. In fact, since $y \sim r^j, dy/dr \sim jr^{j-1}$
as $r\rightarrow 0$, we have $y>0, dy/dr>0$ for $0<r\ll 1$. Suppose there is $r_1>0$ such that
$\frac{d\psi}{dr}(r)>0$ for $0<r<r_1$ and $\frac{d\psi}{dr}(r_1)=0$. Then of course $\psi(r_1)>0$. Moreover, by the equation $(E_j)$, we have
$$\frac{d^2\psi}{dr^2}(r_1) =(j(j+1)-\nu\theta_{\sharp}(r)^{\nu-1}r^2)|_{r=r_1}\frac{\psi(r_1)}{r_1^2}>0,$$
a contradiction. Thus $\frac{d\psi}{dr}(r)>0$ everywhere.
\end{Remark}

Of course, a rigorous result is called for to validate Lemma \ref{Lemma1}. To that end, we will prove the following that justifies A. Kovetz's numerical result as described in Remark \ref{Rem1}:

\begin{Theorem}\label{Thm1}
If $2\leq \nu <5$, then we have
\begin{equation}
\nu\theta(r;\nu)^{\nu-1}r^2 < 6
\end{equation}
for $0\leq r \leq \xi_1(\nu)$.
\end{Theorem}

Theorem \ref{Thm1} leads to an important corollary: 

\begin{Corollary}
Let $2\leq \nu <5$. Then the assertion of Lemma \ref{Lem1}
holds for $j\geq 2$, and Lemma \ref{Lemma1} is verified.
Moreover the solution $y=\psi(r)$ of the equation $(E_j)$ such that
$y\sim r^j$ as $r\rightarrow 0$ is positive and monotone increasing in $r$, a fortiori,
$$\frac{j+1}{r}y+\frac{dy}{dr}>0$$
for $y=\psi(r), r=\xi_1(\nu)$, provided that $j\geq 2$.
\end{Corollary}

\noindent{\it Proof of Theorem \ref{Thm1}.} 
Let us consider
$$g(r):=r^2\theta^{\nu-1}(r)$$
and
$$\bar{m}(\nu):=\sup \nu g(r).$$

Let $r_1\in (0,\xi_1(\nu))$ attain the maximum of $g$. Then $\frac{dg}{dr}(r_1)=0$.
Since
\begin{equation}
\frac{dg}{dr}(r)=r^2\theta^{\nu-1}(r)\Big(
(\nu-1)\frac{\frac{d\theta}{dr}}{\theta}+\frac{2}{r}\Big),
\end{equation}
we have
\begin{equation}
\frac{2}{\nu-1}\theta(r_1)r_1=
-\frac{d\theta}{dr}(r_1)r_1^2=
\int_0^{r_1}\theta^{\nu}(r)r^2dr.\label{L2}
\end{equation}

We use the following fact\\

{ (*): $f(r):=\displaystyle -\frac{r\frac{d\theta}{dr}(r)}{\theta(r)}$ is increasing in $r\in [0,\xi_1(\nu))$, provided that $1<\nu<5$. }\\

This implies $\frac{d}{dr}g(r)\geq 0$ for $0\leq r\leq r_1$. The fact {(*)} can be shown by the phase portrait analysis of the plane dynamical system
$$r\frac{dv}{dr}=-v+v^2+w,\quad r\frac{dw}{dr}=w(2-(\nu-1)v))$$
for $$ v:=-\frac{r}{u}\frac{du}{dr},\quad w:=r^2 u^{\nu-1}$$
of the Lane-Emden equation
$$\frac{d^2u}{dr^2}+\frac{2}{r}\frac{du}{dr}+u^{\nu}=0.$$
Here we are looking $v=f(r), w=g(r)$ along the Lane-Emden orbit
$u=\theta(r;\nu)$.
Note that, if $\nu>5$, $\xi_1=+\infty$ and $f(r)$ oscillates. For the proof, see \cite{JosephL}. 

The following identity will play a key role:
\begin{equation}
\int_0^{r_1}\theta^{\nu}(r)r^2dr=
\frac{r_1^3}{3}\theta^{\nu}(r_1)+
\frac{\nu}{6}\Big(\frac{2}{\nu-1}\Big)^2\theta(r_1)r_1+\frac{\nu}{6}Q,\label{L3}
\end{equation}
where
\begin{equation}
Q:=\int_0^{r_1}\frac{r\frac{d\theta}{dr}(r)+\theta(r)}{\theta^2(r)r^2}
\Big(\int_0^r\theta^{\nu}(s)s^2ds\Big)^2dr. \label{L4}
\end{equation}

Let us show this. 

\begin{align*}
&\int_0^{r_1}\theta^{\nu}(r)r^2dr=
\frac{r_1^3}{3}\theta^{\nu}(r_1)-\frac{\nu}{3}\int_0^{r_1}
\theta^{\nu-1}(r)\frac{d}{dr}\theta(r)r^2dr \\
&=\frac{r_1^3}{3}\theta^{\nu}(r_1)
+\frac{\nu}{3}\int_0^{r_1}\frac{\theta^{\nu}(r)r^2}{\theta(r)r}\Big(
\int_0^{r}\theta^{\nu}(s)s^2ds\Big)dr \\
&=\frac{r_1^3}{3}\theta^{\nu}(r_1)+
\frac{\nu}{6}\int_0^{r_1}
\frac{1}{\theta(r)r}\frac{d}{dr}\Big(\int_0^r\theta^{\nu}(s)s^2ds\Big)^2dr \\
&=\frac{r_1^3}{3}\theta^{\nu}(r_1)
+\frac{\nu}{6}\frac{1}{\theta(r_1)r_1}
\Big(\int_0^{r_1}\theta^{\nu}(s)s^2ds\Big)^2
-\frac{\nu}{6}\int_0^{r_1}
\frac{d}{dr}\Big(\frac{1}{\theta(r)r}\Big)\Big(
\int_0^r\theta^{\nu}(s)s^2ds\Big)^2dr \\
&=\frac{r_1^3}{3}\theta^{\nu}(r_1)
+\frac{\nu}{6}\frac{r_1}{\theta(r_1)}\Big(
\frac{2}{\nu-1}\theta(r_1)\Big)^2+
\frac{\nu}{6}
\int_0^{r_1}
\frac{\frac{d\theta}{dr}(r)r+\theta(r)}{\theta(r)^2r^2}
\Big(\int_0^r\theta^{\nu}(s)s^2ds\Big)^2dr.
\end{align*}
This is \eqref{L3} and \eqref{L4}.

We now divide into two cases. 

1): Suppose $3\leq \nu$.

Then we have
$$\frac{d\theta}{dr}(r)r+\theta(r) \geq \frac{d\theta}{dr}(r)r+\frac{2}{\nu-1}\theta(r)\geq 0,$$
since $\frac{d}{dr}g(r)\geq 0$ for $0\leq r\leq r_1$. So, $Q\geq 0$ and \eqref{L2}, \eqref{L3}
imply 
$$\frac{2}{\nu-1}\geq \frac{g(r_1)}{3}+
\frac{\nu}{6}\Big(\frac{2}{\nu-1}\Big)^2, $$
or
$$g(r_1)\leq \frac{4\nu-6}{(\nu-1)^2}<\frac{6}{\nu}.$$

2): Suppose $1<\nu<3$.

Then,
since $$\frac{d}{dr}\Big(r\frac{d\theta}{dr}+\theta\Big)=
-r\theta^{\nu}<0,$$ we have, for $0\leq r\leq r_1$, \eqref{L2} implies
$$r\frac{d\theta}{dr}(r)+\theta(r)\geq
r_1\frac{d\theta}{dr}(r_1)+\theta(r_1)=
\frac{\nu-3}{\nu-1}\theta(r_1).$$
Therefore we have
\begin{equation}
Q\geq\frac{\nu-3}{\nu-1}\theta(r_1)
\int_0^{r_1}
\frac{1}{\theta^2(r)r^2}
\Big(\int_0^r\theta^{\nu}(s)s^2ds\Big)^2dr.
\end{equation}
On the other hand, 
\begin{align*}
\int_0^{r_1}\frac{1}{\theta^2(r)r^2}
\Big(\int_0^r\theta(s)^{\nu}s^2ds\Big)^2dr &=
\int_0^{r_1}\Big(\frac{r\frac{d}{dr}\theta(r)}{\theta(r)}\Big)^2dr \\
&=-\int_0^{r_1}
r^2\frac{d}{dr}\theta(r)\frac{d}{dr}\Big(\frac{1}{\theta}\Big)(r)dr \\
&=-r_1^2\frac{\frac{d}{dr}\theta(r_1)}{\theta(r_1)}+
\int_0^{r_1}
\frac{1}{\theta(r)}\frac{d}{dr}(r^2\frac{d}{dr}\theta(r))dr \\
&=-r_1^2\frac{\frac{d}{dr}\theta(r_1)}{\theta(r_1)}
-\int_0^{r_1}r^2\theta^{\nu-1}(r)dr \\
&=\frac{2}{\nu-1}r_1-\int_0^{r_1}r^2\theta^{\nu-1}(r)dr\\
&<\frac{2}{\nu-1}r_1-\frac{r_1^3}{3}\theta^{\nu-1}(r_1).
\end{align*}
Therefore
\begin{align*}
Q&>-\frac{3-\nu}{\nu-1}\frac{2}{\nu-1}r_1\theta(r_1)+
\frac{3-\nu}{\nu-1}\frac{r_1^3}{3}\theta^{\nu}(r_1) \\
&=-\frac{3-\nu}{\nu-1}\frac{2}{\nu-1}r_1\theta(r_1)+
\frac{3-\nu}{\nu-1}\frac{r_1^3}{3}\theta(r_1)g(r_1).
\end{align*}
Inserting this estimate to \eqref{L2}, \eqref{L3}, we get
\begin{align*}
\frac{2}{\nu-1}&>\frac{1}{3}g(r_1)+
\frac{\nu}{6}\Big(\frac{2}{\nu-1}\Big)^2-
\frac{\nu}{6}
\frac{3-\nu}{\nu-1}\frac{2}{\nu-1} \\
&+\frac{\nu}{6}\frac{3-\nu}{\nu-1}\frac{1}{3}g(r_1),
\end{align*}
or
$$\nu g(r_1)<M:=\frac{3\nu(7\nu-6-\nu^2)}{(\nu-1)(9\nu-6-\nu^2)}.
$$
For $2\leq\nu$, we see clearly $\displaystyle M<\frac{3\nu}{\nu-1}\leq 6$.
 This completes the proof. $\square$\\

We have shown the following:

\begin{Proposition}\label{Prop5} Let $2\leq \nu <5$. Then
there exists the bounded linear inverse operator
$(1-D\mathcal{G}(\theta))^{-1} \in \mathcal{B}(\mathfrak{E})$.
\end{Proposition}

In this context, hereafter we suppose\\

{\bf (AS1):} \hspace{8mm} $2\leq \nu < 5$, that is,
$\frac{6}{5}<\gamma\leq \frac{3}{2}$.\\

\begin{Remark}
In Theorem \ref{Thm1}, we have not tried to optimize a lower bound of $\nu$ and in fact $\nu$ can go below 2, but not by far from the current argument. Nevertheless, the lower bound in Theorem \ref{Thm1} is not the main  reason of the assumption on $\nu$ in Theorem \ref{Theorem1}. See Section \ref{Sec:dis} for further discussion. 
\end{Remark}

\section{Existence of `distorted Lane-Emden functions': Proof of Theorem \ref{Theorem1}}\label{Sec:5}

We are seeking a solution $u\in\mathfrak{E}$ of the
equation
$$u=\varepsilon\mathfrak{g}+\mathcal{G}(u) \eqno(\ref{6})$$
of the form
\begin{equation}\label{u}
u=\theta +\varepsilon w,
\end{equation}
where $\varepsilon >0$ is sufficiently small and
$\|w\|_{\mathfrak{E}}\leq \Lambda$, $\Lambda$ being specified later.\\

We put 
\begin{equation}\label{g}
\mathcal{G}(\theta+h)=\mathcal{G}(\theta)+
D\mathcal{G}(\theta)h+\omega(h),
\end{equation}
where
\begin{equation}
\omega(h)=\int_0^1(
D\mathcal{G}(\theta+th)-D\mathcal{G}(\theta))h dt. \label{23}
\end{equation}

Recall that $\Lambda_0\geq 2\max\{1, \mu_1(\nu)/\xi_1(\nu)\}$. Similarly to Lemma \ref{Prop3}, we obtain the following: 

\begin{Lemma}\label{Prop6}
If $|u|, |u+h|\leq\Lambda_0$, then 
\begin{equation}
|(u+h)_{\sharp}^{\nu-1}-u_{\sharp}^{\nu-1}|\leq C
|h|^{(\nu-1)\wedge 1}, 
\end{equation}
for some constant $C>0$. 
\end{Lemma}

Thus from the definition of $\omega(h)$ in \eqref{23}, \eqref{defDg} and \eqref{est2} in Proposition \ref{Proposition1}, and Lemma \ref{Prop6} we deduce that

\begin{Proposition}\label{Prop7} 
If $\|h\|_{\mathfrak{E}}\leq \Lambda_0/2$, then
we have
\begin{equation}
\|\omega(h)\|_{\mathfrak{E}}\leq C\|h\|_{\mathfrak{E}}^{\nu\wedge 2}, 
\label{49}
\end{equation}
for some constant $C>0$. 
Moreover, if $\|h_{\ell}\|_{\mathfrak{E}}\leq\Lambda_0/2$ for $
\ell=1,2$, then we have
\begin{align}
\|\omega(h_2)-\omega(h_1)\|_{\mathfrak{E}}&\leq
C\Big(\|h_2-h_1\|_{\mathfrak{E}}^{(\nu-1)\wedge 1}
\|h_2\|_{\mathfrak{E}} + \nonumber \\
&+\|h_1\|_{\mathfrak{E}}^{(\nu-1)\wedge 1}\|h_2-h_1\|_{\mathfrak{E}} \Big), 
\end{align}
for some constant $C>0$. 
\end{Proposition}

By \eqref{u} and \eqref{g}, the equation \eqref{6} is rewritten as 
$$\theta+\varepsilon w
=
\varepsilon\mathfrak{g}+
\mathcal{G}(\theta)+
\varepsilon D\mathcal{G}(\theta)w+\omega(\varepsilon w). $$
Since $\theta=\mathcal{G}(\theta)$, the problem to be solved is reduced to
\begin{equation}
w=\mathfrak{T}(w),
\end{equation}
where
\begin{equation}
\mathfrak{T}(w):=(1-D\mathcal{G}(\theta))^{-1}\Big(
\mathfrak{g}+\frac{1}{\varepsilon}\omega(\varepsilon w)\Big). \label{52}
\end{equation}.\\

Now we propose to solve this problem by the iteration
$$w^{(0)}=0, \qquad w^{(n+1)}=\mathfrak{T}(w^{(n)}).$$
 
In fact let us choose $\Lambda$ so that
$$\|(1-D\mathcal{G}(\theta))^{-1}\mathfrak{g}\|_{\mathfrak{E}}\leq \frac{\Lambda}{2}.$$
By Proposition \ref{Prop7} we see that 
$$\|(1-D\mathcal{G}(\theta))^{-1}\frac{1}{\varepsilon}
\omega(\varepsilon w)\|_{\mathfrak{E}}
\leq C \varepsilon^{(\nu-1)\wedge1}\|w\|_{\mathfrak{E}}\leq
\frac{\Lambda}{2}$$
for $\|w\|_{\mathfrak{E}}\leq\Lambda$ and $\varepsilon \leq \epsilon_1$,
$\epsilon_1$ being sufficiently small. \\

Of course we assume that $\|\varepsilon w\|_{\mathfrak{E}}\leq \delta_0$ for
$\|w\|_{\mathfrak{E}}\leq \Lambda$ and $\varepsilon \leq \epsilon_1$.
Then $\mathfrak{T}$ maps $\mathfrak{X}:=\{ u\in\mathfrak{E} : \|u\|_{\mathfrak{E}}
\leq \Lambda \} $ into itself, provided that $\varepsilon \leq\epsilon_1$.
Moreover, by Proposition \ref{Prop7}, we have 
$$\|\mathfrak{T}(w_2)-\mathfrak{T}(w_1)\|_{\mathfrak{E}}
\leq C \varepsilon^{(\nu-1)\wedge 1}
\|w_2-w_1\|_{\mathfrak{E}}^{(\nu-1)\wedge 1}$$
for $w_1,w_2 \in \mathfrak{X}$ and $\varepsilon \leq\epsilon_1$.

Hence we deduce that $\mathfrak{T}$ is a contraction with respect to the norm $\|\cdot\|_{\mathfrak{E}}$
by taking $\epsilon_1$ smaller if necessary,
under the stronger assumption
{\bf (AS1)}, which guarantees $\nu\geq 2$ so that
$(\nu-1)\wedge 1=1$. 

Therefore, under {\bf (AS1)} we get a fixed point $w$
of $\mathfrak{T}$ in $\mathfrak{X}\subset \mathfrak E$.  Let us denote this unique solution $w$ by $w(r,\zeta;\nu,\varepsilon)$. 

\begin{Definition} Let $$\Theta(r,\zeta;\nu,\varepsilon)
:=\theta(r;\nu)+\varepsilon w(r,\zeta;\nu,\varepsilon).$$ Of course we put $\Theta(r, \zeta;\nu,0)=\theta(r;\nu)$.
We call $\Theta(\cdot,\cdot;\nu,\varepsilon)$  ``the distorted
Lane-Emden function of index $\nu$ with the squared angular
velocity parameter $\varepsilon$". 
\end{Definition}

Note that we can extend $\Theta$ for $r>R_0$. Originally
$w(r,\zeta;\nu,\varepsilon)=(\Theta-\theta)/\varepsilon$ and
$\Theta(r,\zeta;\nu,\varepsilon)$ are defined for $r\leq R_0$.
But $\Theta=\theta+\varepsilon w$ satisfies
\begin{equation}
\Theta=\varepsilon\mathfrak{g}+1+
\mathcal{K}\Theta_{\sharp}^\nu
-(\mathcal{K}\Theta_{\sharp}^\nu)(0,0).\label{30}
\end{equation}
Since $\mathfrak{g}$ and $\mathcal{K}\Theta_{\sharp}$ are
defined on $[0,+\infty)\times[-1,1]$, we can fix the extension of
$\Theta$ by the right-hand side of \eqref{30}. \\

The distorted Lane-Emden function $\Theta$ is the solution of our problem. From the fact that $\|w\|_{\mathfrak{E}}\leq\Lambda$, we easily derive that 

\begin{Proposition}\label{Prop8}
We have
\begin{equation}
|\Theta(r,\zeta;\nu,\varepsilon)-\theta(r;\nu)|\leq C\varepsilon
\end{equation}
for $0\leq r\leq 2\xi_1(\nu), -1\leq\zeta\leq 1, 0<\varepsilon\leq\epsilon(\nu)$.
\end{Proposition}

This finishes the proof of Theorem \ref{Theorem1}.

\section{Properties of the distorted Lane-Emden functions: Proof of Theorem \ref{Theorem2}}\label{Sec:6}

This section is devoted to the properties of our solution $\Theta$. 

\subsection{Regularity of $\Theta$}

\begin{Proposition}\label{Prop9}
We have
\begin{equation}
\Big(\frac{\partial\Theta}{\partial r}-\frac{d\theta}{dr}\Big)^2+
\frac{1-\zeta^2}{r^2}\Big(\frac{\partial\Theta}{\partial\zeta}\Big)^2
\leq C\varepsilon^2
\end{equation}
for $\Theta=\Theta(r,\zeta;\nu,\varepsilon),
\theta=\theta(r;\nu), 0\leq r\leq 2\xi_1(\nu), -1\leq\zeta\leq 1,
0<\varepsilon\leq\epsilon(\nu)$.
\end{Proposition}

\noindent{\it Proof of Proposition \ref{Prop9}.} Recall that $w$ satisfies
\begin{equation}
w=\mathfrak{g}+D\mathcal{G}(\theta)w+
\frac{1}{\varepsilon}\omega(\varepsilon w).
\end{equation}
On the other hand, by \eqref{18}
we have
$$\|\nabla D\mathcal{G}(\theta)w\|_{\mathfrak{E}}=
\|\nabla \mathcal{K}(\nu\theta_{\sharp}^{\nu-1}w)\|_{\mathfrak{E}}
\leq C\|w\|_{\mathfrak{E}}$$
and similarly, 
$$\|\nabla\omega(\varepsilon w)\|_{\mathfrak{E}}\leq C\varepsilon \|w\|_{\mathfrak{E}}
$$
where we recall \eqref{23}. 
Moreover, it is easy to check 
$$\|\nabla\mathfrak{g}\|_{\mathfrak{E}}=\sup\frac{1}{2}\sqrt{1-\zeta^2}r =\frac{R_0}{2}.$$
Therefore we have
$$\|\nabla w\|_{\mathfrak{E}}\leq C, $$
that is, 
\begin{equation}
\|\nabla\Theta -\nabla\theta\|_{\mathfrak{E}}
\leq C\varepsilon. \label{B2}
\end{equation}
$\square$

Here we note that
\begin{align*}
\nabla U&=\left(\frac{\partial U}{\partial x}, \frac{\partial U}{\partial y}, \frac{\partial U}{\partial z}\right)^T \\
&=
\begin{bmatrix}
(\displaystyle \sqrt{1-\zeta^2}\cos\phi \frac{\partial U}{\partial r}
-\frac{\zeta\sqrt{1-\zeta^2}}{r}\frac{\partial U}{\partial \zeta} \\
\\
\displaystyle \sqrt{1-\zeta^2}\sin\phi\frac{\partial U}{\partial r}
-\frac{\zeta\sqrt{1-\zeta^2}}{r}\frac{\partial U}{\partial \zeta} \\
\\
\displaystyle \zeta\frac{\partial U}{\partial r}+\frac{1-\zeta^2}{r}\frac{\partial U}{\partial\zeta}
\end{bmatrix}
\end{align*}
and
$$
|\nabla U|^2=\Big(\frac{\partial U}{\partial r}\Big)^2+\frac{1-\zeta^2}{r^2}\Big(\frac{\partial U}{\partial \zeta}\Big)^2.
$$

Since $\Theta_\sharp^\nu$ is H\"{o}lder continuous, from the standard elliptic theory, it follows that $\Theta \in C^{2,\alpha}$. This establishes (1) of Theorem \ref{Theorem2}.

\subsection{Boundary surface}

Let us fix $r_0$ such that $0<r_0<\xi_1(\nu)$. Then
\begin{equation}
0<\theta_0:=\theta(r_0;\nu)\leq\theta(r;\nu)
\end{equation}
for $0\leq r\leq r_0$. Thanks to Proposition \ref{Prop8}, we can claim, for a sufficiently small $\epsilon_0(\leq\epsilon(\nu))$, that
\begin{equation}
0<\frac{\theta_0}{2}\leq\Theta(r,\zeta;\nu,\varepsilon)
\end{equation}
for $0\leq r\leq r_0, -1\leq \zeta\leq 1, 0,\varepsilon\leq\epsilon_0$.\\

Let us fix $r_1$ such that $\xi_1(\nu) <r_1\leq 2\xi_1(\nu)$. Then
\begin{equation}
\theta(r;\nu)\leq -|\theta_1|<0
\end{equation}
for $r_1\leq r\leq 2\xi_1(\nu)$. 
Here $\theta_1:=\theta(r_1,\nu)$.  Thanks to Proposition \ref{Prop8}, we can claim, taking $\epsilon_0$ smaller if necessary, that
\begin{equation}
\Theta(r,\zeta;\nu,\varepsilon)\leq -\frac{|\theta_1|}{2}<0
\end{equation}
for $r_1\leq r\leq 2\xi_1(\nu), -1\leq \zeta\leq 1, 0<
\varepsilon\leq\epsilon_0$.

Recall that $d\theta/dr <0$ for $0<r<+\infty$. Therefore
there is a positive number $\kappa$ such that
\begin{equation}
\frac{d\theta}{dr}(r;\nu) \leq -\kappa<0
\end{equation}
for $r_0\leq r \leq r_1$.  Thanks to Proposition \ref{Prop9}, we can claim,
taking $\epsilon_0$ smaller if necessary, that
\begin{equation}
\frac{\partial\Theta}{\partial r}(r,\zeta;\nu,\varepsilon)\leq -\frac{\kappa}{2}<0
\end{equation}
for $r_0\leq r\leq r_1, -1\leq\zeta\leq 1, 0<
\varepsilon\leq\epsilon_0$.\\

As a consequence, there exists a unique $\Xi_1(\zeta)=\Xi_1(\zeta;\nu,
\varepsilon)$ in $(r_0,r_1)$ for each $\zeta\in [-1,1]$ so that  
\begin{align*}
&\Theta(r,\zeta;\nu,\varepsilon)>0 \quad\mbox{for}\quad
0\leq r <\Xi_1(\zeta), \nonumber \\
&\Theta(\Xi_1(\zeta),\zeta;\nu,\varepsilon)=0, \nonumber \\
&\Theta(r,\zeta;\nu,\varepsilon)<0\quad\mbox{for}\quad
\Xi_1(\zeta)<r\leq 2\xi_1. \label{B3}
\end{align*}
for $0<\varepsilon\leq\epsilon_0$.

Note that
$$
|\theta(\Xi_1(\zeta))-\theta(\xi_1)|=
|\theta(\Xi_1(\zeta))|=
|\Theta(\Xi_1(\zeta),\zeta)-\theta(\Xi_1(\zeta))|\leq C\varepsilon$$
and $\frac{d\theta}{dr}\leq -\frac{1}{C} <0$ for $r_0\leq r\leq r_1$. Therefore we have
$$|\Xi_1(\zeta) - \xi_1|\leq C\varepsilon.$$

In the same way we can show that the function 
$\zeta\mapsto \Xi_1(\zeta)$ is continuous on $[-1,1]$. In fact, for
$\zeta_{\ell}\in [-1, 1],\ell=1,2$, we look at
$$|\Theta(\Xi_1(\zeta_2),\zeta_2)-\Theta(\Xi_1(\zeta_1),\zeta_2)|=
|\Theta(\Xi_1(\zeta_1),\zeta_2)-
\Theta(\Xi_1(\zeta_1),\zeta_1)|. $$
But 
$$\mbox{(the left-hand side)} \geq \frac{1}{C}|\Xi_1(\zeta_2)-\Xi_1(\zeta_1)|, $$
since 
$$\frac{\partial\Theta}{\partial r}\leq \frac{d\theta}{dr}+C\varepsilon 
\leq-\frac{1}{C'}$$
thanks to Proposition \ref{Prop9}. Since $\zeta\mapsto \Theta(r,\zeta)$ is continuous,
we see 

\noindent $|\Xi_1(\zeta_2)-\Xi_1(\zeta_1)|\rightarrow 0$
as $\zeta_2\rightarrow \zeta_1$.

Moreover the function
$\zeta\mapsto \Xi_1(\zeta)$ is continuously differentiable in 
$(-1,1)$ and the derivative is estimated as
\begin{equation}
\Big|\frac{d\Xi_1}{d\zeta}\Big|=
\Big|-\Big(\frac{\partial\Theta}{\partial r}\Big)^{-1}
\frac{\partial\Theta}{\partial\zeta}\Big|_{r=\Xi_1}\leq C\frac{\varepsilon}{\sqrt{1-\zeta^2}}. \label{47}
\end{equation}

Of course, $\Xi_1(\zeta)=\Xi_1(-\zeta)$. This completes the proof of (2) of Theorem \ref{Theorem2}. 

\subsection{Physical vacuum boundary}

Physical vacuum boundary condition for compressible gas is defined by
\[
-\infty<  \frac{\partial u}{\partial \vec{N}} \Big|_{\partial \mathfrak D}<0
\]
where $u$ is the enthalpy given in \eqref{enthalpy} and $\vec N$ is the unit outer normal vector along $\partial \mathfrak D$ \cite{JaMa2015,Makino2016}.  We will verify this condition for our solution. 
Let us observe the boundary of the domain
\begin{equation}
\mathfrak{D}:=\{(r,\zeta): \Theta(r,\zeta;\nu,\varepsilon)>0\}=
\{(r,\zeta) : H(r,\zeta)<0\},
\end{equation}
where
\begin{equation}
H(r,\zeta):=r-\Xi_1(\zeta;\nu,\varepsilon).
\end{equation}
Let $\mathrm{P}(\Xi_1(\zeta_0),\zeta_0)$ be a boundary point of
$\mathfrak{D}$. The unit outer normal vector 
$\vec{N}$ at $\mathrm{P}$ is given by
$$\vec{N}=\frac{\nabla H}{|\nabla H|},$$
where
\begin{align*}
&\nabla H=
\begin{bmatrix}
\displaystyle \frac{x}{r}+\frac{x\zeta}{r^2}\frac{d\Xi_1}{d\zeta} \\
\\
\displaystyle \frac{y}{r}+\frac{y\zeta}{r^2}\frac{d\Xi_1}{d\zeta} \\
\\
\displaystyle \zeta-\frac{1-\zeta^2}{r}\frac{d\Xi_1}{d\zeta}
\end{bmatrix}, \\
&|\nabla H|^2=1+\frac{1-\zeta^2}{r^2}\Big(\frac{d\Xi_1}{d\zeta}\Big)^2.
\end{align*}
Note that $|\nabla H|=1+O(\varepsilon^2)$ thanks to \eqref{47}.
Moreover we have
\begin{align*}
\frac{\partial \Theta}{\partial\vec{N}}&=
\frac{\displaystyle (\nabla H|\nabla\Theta)}{|\nabla H|}=
\frac{\displaystyle \frac{\partial\Theta}{\partial r}-
\frac{1-\zeta^2}{r^2}\frac{d\Xi_1}{d\zeta}\frac{\partial\Theta}{\partial\zeta}}
{\displaystyle\sqrt{ 1+\frac{1-\zeta^2}{r^2}\Big(\frac{d\Xi_1}{d\zeta}\Big)^2}} = \\
&=\frac{\partial\Theta}{\partial r}\Big|_{\mathrm{P}}+O(\varepsilon^2)
=\frac{d\theta}{dr}\Big|_{r=\Xi_1(\zeta_0)}+O(\varepsilon) = \\
&=-\frac{\mu_1(\nu)}{\xi_1(\nu)^2}+O(\varepsilon).
\end{align*}
In this sense, taking $\epsilon_0$ smaller if necessary, we can claim that 
{the boundary point $\mathrm{P}$ is a physical vacuum boundary}, provided that $0<\varepsilon\leq\epsilon_0$.

\subsection{`Old caricature' presented by S. Chandrasekhar}

We are considering the domain
\begin{equation}
\mathfrak{D}:=\{\vec{x}\in\mathbb{R}^3 : r=|\vec{x}|\leq R_0, \Theta(r,\zeta)>0\},
\end{equation}
which is occupied by the rotating gaseous star, and its surface
\begin{equation}
\partial\mathfrak{D}=\{\vec{x}\in\mathbb{R}^3 : |\vec{x}|\leq R_0,
\Theta(r,\zeta)=0 \}.
\end{equation}
We can claim that $\partial\mathfrak{D}$ is a $C^1$-submanifold of $\mathbb{R}^3$, since 
$\Theta$ is, as a function of $\vec{x}=(x,y,z)$, of class $C^1$ on $\mathbb{R}^3$ and $\nabla\Theta$ does not vanish on the surface. 
 In fact, $\Theta$ satisfies
$$\|\nabla\Theta-\nabla\theta\|_{\mathfrak{E}}\leq C\varepsilon,
$$
and since 
$$\nabla\theta =\frac{\vec{x}}{r}\frac{d\theta}{dr} \not=\vec{0}
$$ 
we have $$\frac{1}{C}\leq \|\nabla\Theta\|_{L^{\infty}}\leq C
$$
for $r_0\leq r\leq r_1, 0<\varepsilon\leq\epsilon_1$. Therefore the old caricature quoted in \cite[p. 253]{Chandrasekhar1967} as `Figure 2'
is misleading. In this caricature, the North Pole and the South Pole look being peaked. It is impossible, since $\partial\mathfrak{D}$ is axisymmetric and of class $C^1$. For the sake of certainty, we observe
the configuration of the curve $\mathfrak{S}$ of 
$\partial\mathfrak{D}=\{r=\Xi_1(\zeta)\}$ in the $(\varpi, z)$-plane.

Consider $\mathfrak{S}$ near the North Pole $\mathrm{P}:
(\varpi,z)=(0, \Xi_1(1))$. Then it can be described by a function
\begin{equation}
z=Z(\varpi), \qquad 0\leq \varpi\ll 1.
\end{equation}
In fact the equation 
$$\varpi=\Xi_1(\zeta)\sqrt{1-\zeta^2}$$
admits a continuously differentiable inverse function
$\zeta=\zeta(\varpi)$ for $0<\varpi\ll 1$ near $\mathrm{P}$.
In fact we see $d\varpi/d\zeta\rightarrow -\infty$ as $\zeta
\rightarrow 1-0$. Here we have used 
$$\sqrt{1-\zeta^2}\frac{d\Xi_1}{d\zeta}=-\Big(\frac{\partial\Theta}{\partial r}\Big)^{-1}\sqrt{1-\zeta^2}
\frac{\partial\Theta}{\partial \zeta}=O(\varepsilon).$$ By this inverse
$\zeta=\zeta(\varpi)$, we put
$$Z(\varpi)=\Xi_1(\zeta(\varpi))\zeta(\varpi). $$
Then, on the contrary to the impression by the old caricature 
of \cite[Figure 2]{Chandrasekhar1967}, we claim

\begin{Proposition}\label{Prop11}
We have
\begin{equation}
\frac{dZ}{d\varpi}\rightarrow 0\qquad\mbox{as}\quad \varpi\rightarrow 0.
\end{equation}
\end{Proposition}

In order to verify Proposition \ref{Prop11}, we prepare the following

\begin{Proposition}\label{Prop12}
For any $r_*>0$ we have
\begin{equation}
\sqrt{1-\zeta^2}\frac{\partial\Theta}{\partial\zeta}(r,\zeta)\rightarrow 0
\quad\mbox{as}\quad
r\rightarrow r_*, \;\; \zeta\rightarrow 1-0.
\end{equation}
\end{Proposition}

Proof. Recall
\begin{align*}
\frac{\partial\Theta}{\partial x}&=\frac{x}{r}\frac{\partial\Theta}{\partial r}-\frac{x\zeta}{r^2}\frac{\partial\Theta}{\partial\zeta} \\
&=\sqrt{1-\zeta^2}(\cos\phi)\frac{\partial\Theta}{\partial r}
-\frac{\zeta}{r}(\cos\phi)\sqrt{1-\zeta^2}
\frac{\partial\Theta}{\partial\zeta}, \\
\frac{\partial\Theta}{\partial y}&=\frac{y}{r}\frac{\partial\Theta}{\partial r}
-\frac{y\zeta}{r^2}\frac{\partial\Theta}{\partial\zeta} \\
&=\sqrt{1-\zeta^2}(\sin\phi)\frac{\partial\Theta}{\partial r}-
\frac{\zeta}{r}(\sin\phi)\sqrt{1-\zeta^2}
\frac{\partial\Theta}{\partial\zeta},
\end{align*}
while $$x=r\sqrt{1-\zeta^2}\cos\phi, \quad
y-r\sqrt{1-\zeta^2}\sin\phi, \quad z=r\zeta. $$
Recall that $\Theta$ is a $C^1$-function of $\vec{x}=(x,y,z)$.
Consider the path
$$x=r\sqrt{1-\zeta^2}\cos\phi \rightarrow 0, y=
r\sqrt{1-\zeta^2}\sin\phi\rightarrow0, z=r\zeta\rightarrow r_*,$$
while $r\rightarrow r_*, \zeta\rightarrow 1-0$, keeping $\phi$ being
constant. From the case $\phi=0$, we see that the existence of
$\lim \partial\Theta/\partial x$ implies the existence of
$$A:=\lim_{r\rightarrow r_*,\zeta\rightarrow 1-0}
 \sqrt{1-\zeta^2}\frac{\partial\Theta}{\partial\zeta}(r,\zeta).
$$
Considering arbitrary $\phi$, we have
\begin{align*}
\frac{\partial\Theta}{\partial x}\Big|_{x=y=0,z=r_*}&=-\frac{1}{r_*}
(\cos\phi) A, \\
\frac{\partial\Theta}{\partial y}\Big|_{x=y=0,z=r_*}&=-\frac{1}{r_*}
(\sin\phi)A.
\end{align*}
Since the left-hand sides are independent of $\phi$, it should be the case
that $A=0$. $\square$\\

\noindent{\it Proof of Proposition \ref{Prop11}.} Look at
\begin{align}
\frac{dZ}{d\varpi}&=
\frac{\zeta\frac{d\Xi_1}{d\zeta}+\Xi_1}{\sqrt{1-\zeta^2}\frac{d\Xi_1}{d\zeta}-\frac{\zeta}{\sqrt{1-\zeta^2}}\Xi_1}\Big|_{\zeta=\zeta(\varpi)} \nonumber \\
& \nonumber\\
&=\frac{\zeta\sqrt{1-\zeta^2}\frac{d\Xi_1}{d\zeta}+\sqrt{1-\zeta^2}\Xi_1}{(1-\zeta^2)\frac{d\Xi_1}{d\zeta}-\zeta\Xi_1}\Big|_{\zeta=\zeta(\varpi)}.
\label{55}
\end{align}
Proposition \ref{Prop12} implies
$$\sqrt{1-\zeta^2}\frac{d\Xi_1}{d\zeta}=-\Big(\frac{\partial\Theta}{\partial r}\Big)^{-1}
\sqrt{1-\zeta^2}\frac{\partial\Theta}{\partial\zeta}\Big|_{r-\Xi_1(\zeta)}
\rightarrow 0$$
as $\zeta\rightarrow 1-0$. The denominator of \eqref{55} tends to
$-\Xi_1(1) \doteqdot -\xi_1 \not= 0$. Thus we get the conclusion.
$\square$

\section{Justification of Chandrasekhar's approximation}\label{sec:justify}

Recall that $w=w(r,\zeta;\nu,\varepsilon)$ is the fixed point
function of the mapping $\mathfrak{T}$ defined by \eqref{52}. Thanks to
\eqref{49}, we see
$$w(r,\zeta;\nu,\varepsilon)=(1-D\mathcal{G}(\theta))^{-1}
\mathfrak{g}+O(\varepsilon),$$
so that
\begin{equation}
\Theta(r,\zeta;\nu,\varepsilon)=\theta(r;\nu)+
\mathfrak{h}(r,\zeta;\nu)\varepsilon+
O(\varepsilon^2), \label{B4}
\end{equation}
where
\begin{equation}
\mathfrak{h}:=(1-D\mathcal{G}(\theta))^{-1}\mathfrak{g}.
\end{equation}

Let us observe the structure of the first approximation $\mathfrak{h}$.

We have
\begin{equation}
-\triangle \mathfrak{h}=\nu\theta_{\sharp}^{\nu-1}\mathfrak{h}-1,\qquad \mathfrak{h}(\vec{0})=0,
\end{equation}
since $\triangle\mathfrak{g}=1$.

Let us consider the expansion by the Legendre's polynomials
\begin{equation}
\mathfrak{h}(r,\zeta)=\sum_{k=0}^{\infty}
h_{2k}(r)P_{2k}(\zeta).
\end{equation}

First, $y=h_0(r)$ is the solution of the equation
\begin{equation}
-\frac{1}{r^2}\frac{d}{dr}r^2\frac{dy}{dr}=\nu\theta_{\sharp}^{\nu-1}y-1
\label{C1}
\end{equation}
such that $y=0$ at $r=0$. However general solution of \eqref{C1} is given by
$$y=y_0(r)+C_1y_1(r)+C_2y_2(r),$$
where 
$y_0(r)$ is a special solution such that $$y_0=\frac{1}{\nu}+[r^2]_1
\quad\mbox{as}\quad r\rightarrow 0,$$
and $y_1(r), y_2(r)$ are solutions of the homogeneous equation such that
$$y_1=1+[r^2]_1,\qquad y_2=\frac{1}{r}(1+[r^2]_1)$$
as $r\rightarrow 0$. Therefore 
$$h_0(r)=y_0(r)-\frac{1}{\nu}y_1(r)$$
is uniquely determined and $h_0(r)=[r^2]_1$.

Next we consider $h_2(r)$. Then $y=h_2(r)$ is a solution 
of the equation
$$\Big[-\frac{1}{r^2}\frac{d}{dr}r^2\frac{d}{dr}+
\frac{j(j+1)}{r^2}\Big]y=\nu\theta_{\sharp}(r)^{\nu-1}y
\leqno(E_j)
$$
for $j=2$ such that $y=0$ at $r=0$. Let $y=\psi_2(r)$ be the solution of $(E_2)$
such that $y=r^2(1+[r^2]_1)$. Then there should exist a constant $A_2$
such that $h_2(r)=A_2\psi_2(r)$, since other independent solution is
$\sim r^{-3}$ as $r\rightarrow 0$. On the other hand $y=h_2(r)$
satisfies
$$ \Big[-\frac{1}{r^2}\frac{d}{dr}r^2\frac{d}{dr}+
\frac{j(j+1)}{r}\Big]y=0 \leqno(H_j) $$
with $j=2$ on $r\geq \xi_1(\nu)$. 
Since
$$\mathfrak{h}-\mathfrak{g}=
D\mathcal{G}(\theta)\mathfrak{h}=\mbox{Const.}+O(r^{-1})$$
as $r\rightarrow+\infty$, and since
$$\mathfrak{g}(r,\zeta)=\frac{r^2}{6}-\frac{r^2}{6}P_2(\zeta),$$
there should exist a constant $C_2$ such taht
 $$h_2(r)+\frac{r^2}{6}=C_2 r^{-3}$$
for $r\geq \xi_1(\nu)$.
Since $h_2(r)$ is continuously differentiable at $r=\xi_1(\nu)$, we have
\begin{subequations}
\begin{align}
-\frac{r^2}{6}+\frac{C_2}{r^3}&=A_2\psi_2(r) \\
-\frac{r}{3}-\frac{3C_2}{r^4}&=A_2D\psi_2(r)
\end{align}
\end{subequations} 
at $r=\xi_1(\nu)$. But we know that
$$\frac{j+1}{r}y+\frac{dy}{dr}>0$$
for $y=\psi_2(r), j=2, r=\xi_1(\nu)$, provided that $2\leq \nu <5$.
Therefore the constant $A_2$ is determined as
\begin{equation}
A_2=-\frac{5}{6}
\frac{r^2}{3\psi_2(r)+r\frac{d}{dr}\psi_2(r)}\Big|_{r=\xi_1(\nu)}.
\end{equation}
Note that $A_2<0$.

Finally we consider $j=2k\geq 4$. Let $y=\psi_j(r)$ be the solution of the
equation $(E_j)$ such that $y\sim r^j$ as $r\rightarrow 0$. Since other independent solution is $\sim r^{-j-1}$, there should exist a constant $A_j$ such that $h_j(r)=A_j\psi_j(r)$. On the other hand 
$h_2(r)$ satisfies $(H_j)$ on $r\geq\xi_1(\nu)$, there should exist
a constant $C_j$ such that $h_j(r)=C_jr^{-j-1}$ on $r\geq \xi_1(\nu)$. Since $h_j$ is continuously
differentiable at $r=\xi_1(\nu)$, we have
\begin{subequations}
\begin{align}
\frac{C_j}{r^{j+1}}&=A_j\psi_j(r) \label{81a}\\
-(j+1)\frac{C_j}{r^{j+2}}&=Aj\frac{d}{dr}\psi_j(r) \label{81b}
\end{align}
\end{subequations}
at $r=\xi_1(\nu)$. But we know 
$$\frac{j+1}{r}y+\frac{dy}{dr}>0$$
for $y=\psi_j(r), j=2k\geq 4, r=\xi_1(\nu)$. Thus the system of linear equations \eqref{81a} \eqref{81b} determines $A_j=C_j=0$
so that $h_j=0$.\\

Summing up, we get
\begin{equation}
\mathfrak{h}(r,\zeta)=h_0(r)+A_2\psi_2(r)P_2(\zeta).
\end{equation}
Recall that $A_2<0$, $\psi_2(r)>0$, and $$P_2(\zeta)=
\frac{1}{2}(3\zeta^2-1).$$

On the other hand, it follows from \eqref{B2}, \eqref{B3}, \eqref{B4} that
\begin{equation}
\Xi_1(\zeta;\varepsilon)=\xi_1+
\frac{\xi_1^2}{\mu_1}\mathfrak{h}(\xi_1,\zeta)\varepsilon 
+O(\varepsilon^2),
\end{equation}
where $\xi_1=\xi_1(\nu), \mu_1=\mu_1(\nu)$ and we note 
${d\theta}/{dr}(\xi_1)=-\mu_1/\xi_1^2$.
Therefore the `oblateness of the surface' defined by
$$\sigma:=\frac{\Xi_1(0;\varepsilon)-\Xi_1(\pm1;\varepsilon)}{\xi_1}$$
turns out to be
$$\sigma=-\frac{3}{2}
\frac{\xi_1}{\mu_1}A_2\psi_2(\xi_1)\varepsilon +
O(\varepsilon^2)>0,
$$
provided that $0<\varepsilon \ll 1$.
In this sense, as for the old caricature in
\cite[p.253]{Chandrasekhar1967}, Newton is better than Cassini.\\

This is the mathematically rigorous justification of
the Chandrasekhar's approximation given in
\cite{Chandrasekhar1933}.

\section{Discussion}\label{Sec:dis}

We have the restriction of $\gamma$ in the interval $(6/5, 3/2]$. This range covers the diatomic gas ($\gamma=7/5$) and the radiation case ($\gamma=4/3$), but excludes the monoatomic gas ($\gamma=5/3$ that is, $\nu=3/2$). Even if we can improve the range of $\gamma$ in 
Theorem \ref{Thm1}, this case is excluded in order that the mapping $\mathfrak{T}$ be a contraction
with respect to the norm $\|\cdot\|_{\mathfrak{E}}$, which requires $\nu \geq 2$. However, we expect the results for the values of $\gamma$ in a wider range  
and the development of other methods will be expected. \\

\

\noindent{\bf\Large Acknowledgment.} This joint work was initiated by the discussion between authors during the stay of the second author at the Korea Institute for Advanced Study on July 11 - 15, 2016. Sincere thanks are addressed to the Korea Institute for Advanced Study for their support and hospitality during this stay. The first author acknowledges the support by NSF grants DMS-1608492 and DMS-1608494. 

\

\end{document}